\documentclass[letterpaper, paper,11pt]{AAS}	% for preprint proceedings

\usepackage{bm}
\usepackage{amsmath}
\usepackage{subfigure}
\usepackage{placeins}
\usepackage[colorlinks=true, pdfstartview=FitV, linkcolor=black, citecolor= black, urlcolor= black]{hyperref}
\usepackage{overcite}
\usepackage{footnpag}			      	% make footnote symbols restart on each page
\usepackage{authblk}
\usepackage[font={small,it}]{caption}

\PaperNumber{21-776}

\begin{document}

\title{Sub-Optimal Fast Fourier Series Approximation for Initial Trajectory Design}

\author{Caleb Gunsaulus\thanks{Undergraduate Student, Department of Aerospace Engineering, Iowa State University, USA.}, Carl De Vries\thanks{Graduate Student, Department of Aerospace Engineering, Iowa State University, USA.}, William Brown\thanks{Undergraduate Student, Department of Aerospace Engineering, Iowa State University, USA.}, Youngro Lee\thanks{Graduate Student, Department of Aerospace Engineering, Iowa State University, USA.}, Madhusudan Vijayakumar\thanks{Graduate Student, Department of Aerospace Engineering, Iowa State University, USA.}, Ossama Abdelkhalik\thanks{Associate Professor, Department of Aerospace Engineering, Iowa State University, Ames, Iowa 50011, AIAA Senior Member.}}
\maketitle{} 		

\begin{abstract}
The Finite Fourier Series (FFS) Shape-Based (SB) trajectory approximation method has been used to rapidly generate initial trajectories that satisfy the dynamics, trajectory boundary conditions, and limitation on maximum thrust acceleration. The FFS SB approach solves a nonlinear programming problem (NLP) in searching for feasible trajectories.
This paper extends the development of the FFS SB approach to generate sub optimal solutions. Specifically, the objective function of the NLP problem is modified to include also a measure for the time of flight. Numerical results presented in this paper show several solutions that differ from those of the original FFS SB ones. The sub-optimal trajectories generated using a time of flight minimization are shown to be physically feasible trajectories and potential candidates for direct solvers.
\end{abstract}

\section{Introduction}
Current trajectory optimization methods can minimize objectives such as time of flight or fuel consumption. Two common methods to solve these problems include direct and indirect methods \cite{doi:10.2514/2.4231}. Direct methods are very useful, but can be computationally expensive, and require initial guesses for the design parameters to start the solution process. Sub-optimal, yet feasible, trajectories can be a good option to initialize a direct solver. However, many sub-optimal methods are useful only in limited scenarios, are subject to geometric requirements, or are not able to satisfy spacecraft thrust constraints. One approximate method is the Finite Fourier Series Shape-Based method \cite{doi:10.2514/1.58789, doi:10.2514/1.G000878, doi:10.1016/j.asr.2015.11.034}. The FFS SB method generates a spacecraft trajectory by solving a NLP problem for a set of Fourier coefficients which approximate the spacecraft states as a function of time. Depending on the problem, some of the Fourier coefficients can be determined from the boundary conditions. The rest of the unknown coefficients are found using an optimizer. In this case, MATLAB's \textit{fmincon} function is used to find trajectories by minimizing an objective function. To find these trajectories, the optimizer computes the error from the FFS SB approximation at multiple points along the trajectory and then minimizes the sum of the squared error. The FFS SB method has been developed for two-body and three-body systems as well as two dimensional and three dimensional spaces \cite{doi:10.2514/1.G000878, doi:10.1016/j.asr.2015.11.034}. The method has been applied to interplanetary transfers, spacecraft rendezvous, orbit raising, and phasing maneuvers. While the FFS SB method doesn’t produce optimal solutions, it was shown to produce cost values (fuel mass) near the optimal solutions obtained using GPOPS in some test cases \cite{doi:10.1016/j.asr.2015.11.034}.

Currently, the FFS SB method does not attempt to minimize objectives such as time of flight or fuel mass. In this paper, a time of flight penalty is added to the FFS SB method objective function. This work solves the multi-objective problem given that the time of flight objective could be competing with the objective of satisfying equations of motion. The approach is applied to spacecraft orbit raising and rendezvous test cases which demonstrate feasible trajectories with shorter flight times are generated using the FFS SB method.

\section{Finite Fourier Series Shape-Based Method}
Mathematically, any given periodic function can be represented as an infinite summation of sine and cosine terms. Shape-Based methods belong to a domain of space trajectory optimization which geometrically approximates the shape of a trajectory. This trajectory is then used to compute the thrust vector and hence the trajectory solution of the design/optimization problem. One such popular and robust representation of the trajectory shape can be obtained by using the Finite Fourier Series. The FFS SB method approximates the states of the spacecraft namely the radius $r$ and polar angle $\theta$ as follows:

\begin{equation}
\label{eq:FFS_r}
    r(t) = \frac{a_0}{2} + \sum_{n=1}^{n_r} [a_n cos( \frac{n\pi t}{T}) + b_n sin( \frac{n\pi t}{T}) ]
 \end{equation}

\begin{equation}
\label{eq:FFS_theta}
     \theta(t) = \frac{c_0}{2} + \sum_{n=1}^{n_\theta} [c_n cos( \frac{n\pi t}{T}) + d_n sin( \frac{n\pi t}{T}) ]
\end{equation}

where, $n_r$ and $n_\theta$ are the number of Fourier coefficients used to represent the radius $r$ and polar angle $\theta$ of the spacecraft states, respectively. The biggest advantage of the FFS method is that for every different selection of the Fourier coefficients, a different trajectory shape is obtained. Additionally, one should note that the states of the trajectory are discretized at 'n' points and are only approximated using the FFS at these discrete points. If the initial and final boundary conditions (BCs) of the spacecraft states are available, they can be used to solve for some of the Fourier coefficients in terms of the others.

%     \begin{align*}
%     a_1 &= \frac{r_i-r_f}{2} - \sum_{n=3}^{n_r} a_n; n_r \geq 3, n:odd &&
%     a_2 &= \frac{r_i+r_f-a_0}{2}\sum_{n=4}^{n_r} a_n; n_r \geq 4, n:even \\
%     b_1 &= \frac{T}{2\pi}(\dot{r_i}-\dot{r_f})-\sum_{n=3}^{n_r} nb_n; n_r \geq 3, n:odd  &&
%     b_2 &= \frac{T}{4\pi}(\dot{r_i}+\dot{r_f})-\sum_{n=4}^{n_r} nb_n; n_r \geq 4, n:even \\
%     c_1 &= \frac{\theta_i-\theta_f}{2} - \sum_{n=3}^{n_\theta} c_n; n_\theta \geq 3, n:odd  &&
%     c_2 &= \frac{\theta_i+\theta_f-c_0}{2}\sum_{n=4}^{n_\theta} c_n; n_\theta \geq 4, n:even\\
%     d_1 &= \frac{T}{2\pi}(\dot{\theta_i}-\dot{\theta_f})-\sum_{n=3}^{n_\theta} nd_n; n_\theta \geq 3, n:odd  &&
%     d_2 &= \frac{T}{4\pi}(\dot{\theta_i}+\dot{\theta_f})-\sum_{n=4}^{n_\theta} nd_n; n_\theta \geq 4, n:even \\
% \end{align*}

The number of Fourier terms to be used is problem dependent and is obtained on a trial and error basis. Once the number of FFS terms and the number of discretization points are determined, the sine and cosine terms become fixed and can be calculated once and used repeatedly to decrease computational cost. The Fourier approximated function can thus be reduced to the following form:
\begin{equation}
    r(t) = [A]*X + B
\end{equation}
where A is a square matrix of order $(2*n_r +1) \times (2*n_r +1)$ consisting of the sine and cosine terms, X is a column matrix having  ($2*n_r +1$) rows representing the unknown Fourier coefficients that have to be solved and B is the column matrix having ($2*n_r +1$) rows that represents the coefficients that were solved for using the BCs.

\section{Equations of Motion}

The polar equations of motions are developed in two dimensions where $r$ is the spacecraft radius and $\theta$ is the polar angle \cite{doi:10.2514/1.36848}. The spacecraft is acted on by a gravitational force from a single, central body and an external thrust input $T_a$. A local vertical local horizontal coordinate frame is fixed to the spacecraft. The angle $\alpha$ is the steering angle measured from the local horizontal and is the angle at which the thrust acceleration is applied. The gravitational parameter is given by $\mu$, and the equations of motion are shown in Eqs. \eqref{eq:r_eom} and \eqref{eq:theta_eom}. 

\begin{equation}
    \label{eq:r_eom}
    \ddot{r}-r\dot{\theta}^2+\frac{\mu}{r^2}=T_asin(\alpha)
\end{equation}
\begin{equation}
    \label{eq:theta_eom}
    2\dot{r}\dot{\theta}+r\ddot{\theta}=T_acos(\alpha)
\end{equation}

The steering angle is commonly equal to the flight path angle, $\gamma$, plus a multiple of $\pi$ because it's allowed to fire in both the positive and negative directions in Eq. \eqref{eq:steering_fpa}  \cite{doi:10.2514/1.58789,doi:10.2514/1.36848,doi:10.2514/8.7261,doi:10.2514/3.56583,doi:10.2514/1.13095,petrothesis}. The additional term $\pi n$ is added to give the direction. The thrust acceleration is applied in the direction of the flight path angle when $n=0$ and opposite to the flight path angle when $n=1$. The steering angle can also be found from the common flight path angle equation as the ratio of the velocity components in the local vertical local horizontal frame shown in Eq. \eqref{eq:fpa}.

\begin{equation}
    \label{eq:steering_fpa}
    \alpha = \gamma + \pi n
\end{equation}

\begin{equation}
    \label{eq:fpa}
    tan(\alpha) = tan(\gamma) = \frac{\dot{r}}{r\dot{\theta}}
\end{equation}

The thrust acceleration shown in Eq. \eqref{eq:ta} can be solved for from Eq. \eqref{eq:theta_eom}. Substituting it into Eq. \eqref{eq:r_eom} and applying Eq. \eqref{eq:fpa} results in a single equation of motion shown in Eq. \eqref{eq:combined_eom}.

\begin{equation}
    \label{eq:ta}
    T_a=\frac{2\dot{r}\dot{\theta} + r\ddot{\theta}}{cos(\alpha)}
\end{equation}
\begin{equation}
    \label{eq:combined_eom}
    F = r^2(\dot{\theta}\ddot{r}-\dot{r}\ddot{\theta}) + \dot{\theta}(\mu - 2r\dot{r}^2)-(r\dot{\theta})^3 = 0
\end{equation}

The FFS SB method exploits the fact that the equations of motion can be formulated in a single expression and be collected easily so all the terms are set equal to zero. The first and second derivatives of the FFS SB state approximations can be computed from Eqs. \eqref{eq:FFS_r} and \eqref{eq:FFS_theta}. When the states and state derivatives are evaluated using specific Fourier coefficients, the resulting values can be plugged into Eq. \eqref{eq:combined_eom}. The result of evaluating this equation will be non-zero because they are approximations. This error quantifies how well the FFS SB approximated values fit the equations of motion. This calculation is performed at multiple points along the trajectory with the goal to minimize the summation of the squared residual error.

\section{Sub-Optimal Solution Formulation}
The objective function developed to solve the multi-objective function is composed of a weighted sum of individual objectives. Two objectives are used and they are 1) the squared error of the FFS approximation to the equations of motion and 2) the time of flight penalty. The objective function is shown in Eq.~\eqref{eq:J_TOF}

\begin{equation}
\label{eq:J_TOF}
    J = (1-\omega) \textbf{F}^T \textbf{F} + \omega \left(\frac{ToF}{ToF_0}\right)
\end{equation}

where, $J$ is the objective function, $\textbf{F}$ is the error in the Fourier series approximation of the equations of motion evaluated multiple points on the trajectory, $ToF$ is the time of flight, $ToF_0$ is the initial time of flight guess, and $\omega$ is the weight applied to the ToF.

Now that the time of flight is included in the objective function, the goal is to select a weight which will minimize the combination of the weighted equation of motion squared error and the weighted time of flight. Modifications must be made to ensure the trajectory's boundary conditions are satisfied since the time of flight is now a design parameter. The final polar angle is a free parameter in the orbit raising formulation, however, spacecraft rendezvous requires a fixed polar angle. The final polar angle at which rendezvous occurs will change throughout the optimization process as the time of flight varies.

The steps to compute the new rendezvous polar angle are shown in Eqs. \eqref{eq:orbital_rate}, \eqref{eq:mean_anomaly}, and \eqref{eq:rendezvous_polar_angle}. The orbital rate for the target orbit can be computed from the gravitational parameter and the semi-major axis for the final, circular orbit. Then the initial guess for the time of flight is subtracted from the current ToF trial value. The change in the time of flight is multiplied by the orbital rate and yields the change in mean anomaly (polar angle) for the target orbit across that time period. Adding the change in mean anomaly to a reference mean anomaly which corresponds to the initial guess for the time of flight results in the new rendezvous polar angle.

\begin{equation}
    \label{eq:orbital_rate}
    n=\sqrt{\frac{\mu}{r_f^3}}
\end{equation}

\begin{equation}
    \label{eq:mean_anomaly}
    \Delta M=n(ToF - ToF_0)
\end{equation}

\begin{equation}
    \label{eq:rendezvous_polar_angle}
    \theta_f = \theta_{f,0} + \Delta M
\end{equation}

\section{Optimal Solution Formulation}

% ---------- GPOPS ----------
% Youngro take here
% ---------------------------
The FFS generates the sub-optimal thrust acceleration by penalizing ToF in the objective function, as shown in the previous design process and results.
In order to examine the similarity between the sub-optimal solutions and optimal solutions, General Purpose OPtimal control Software (GPOPS-II), which is a MATLAB software intended to solve general nonlinear optimal control problems \cite{patterson2014gpops}, is employed to find optimal solutions.
Since GPOPS-II is a direct method to solve the optimal control problem, the equations of motion in Eqs. \eqref{eq:r_eom} and \eqref{eq:theta_eom} are used instead of using the combined equation, which is Eq. \eqref{eq:combined_eom}.
Accordingly, it does not need to consider the motion error itself, so the objective function only includes the ToF term.

The minimum time orbit raising problem can be constructed as shown in Eq. \eqref{eqn:gpops_formulation_OR}.
\begin{subequations}
\label{eqn:gpops_formulation_OR}
\begin{align}
    J & = ToF \\
    \text{subject to   } & Eqs. \eqref{eq:r_eom} \text{ and } \eqref{eq:theta_eom} \\
    & r(0) = r_i, \theta(0) = \theta_i, \dot{r}(0) = \dot{r}_i, \dot{\theta}(0) = \dot{\theta}_i \\
    & r(ToF) = r_f, \dot{r}(ToF) = \dot{r}_f, \dot{\theta}(ToF) = \dot{\theta}_f \\
    & -T_{a, max} \leq T_a \leq T_{a, max} \\
    & r_i \leq r \leq r_f, 0 \leq \theta \leq \theta_{f,0}, 0 \leq \dot{r} \leq 1, 0  \leq \dot{\theta} \leq 1
\end{align}
\end{subequations}
where $\theta_{f,0}$, which is the initial guess of the final polar angle, is set as $12\pi$, and the initial guess of ToF is given as 120000 seconds.
Furthermore, the boundaries of the time derivatives of the radial distance and the polar angle are defined as proper values that provide enough space to compute the optimal trajectory.

GPOPS computes the optimal solution for the rendezvous case as well.
Unlike the orbit raising, the rendezvous case constrains the final polar angle, which can be calculated based on equations \eqref{eq:orbital_rate}, \eqref{eq:mean_anomaly}, and \eqref{eq:rendezvous_polar_angle}. Therefore, the new rendezvous polar angle is included in Eq. \eqref{eqn:gpops_formulation_rdz}.
\begin{subequations}
\label{eqn:gpops_formulation_rdz}
\begin{align}
    J & = ToF \\
    \text{subject to   } & Eqs. \eqref{eq:r_eom} \text{ and } \eqref{eq:theta_eom} \\
    & r(0) = r_i, \theta(0) = \theta_i, \dot{r}(0) = \dot{r}_i, \dot{\theta}(0) = \dot{\theta}_i \\
    & r(ToF) = r_f, \dot{r}(ToF) = \dot{r}_f, \dot{\theta}(ToF) = \dot{\theta}_f \\
    & \theta(ToF) = \theta_{f,0} + \sqrt{\frac{\mu}{r_f^3}} \left( ToF - ToF_0 \right) \\
    & -T_{a, max} \leq T_a \leq T_{a, max} \\
    & r_i \leq r \leq r_f, 0 \leq \theta \leq \theta_{f,0}, 0 \leq \dot{r} \leq 1, 0  \leq \dot{\theta} \leq 1.2
\end{align}
\end{subequations}
The variables $\theta_{f,0}$ and $ToF_0$ are the initial guess for the final polar angle and the time of flight, and their values are given as $13\pi$ and 120000 seconds, respectively.

\section{Test Cases}\label{SecResults}
\subsection{Orbit Raising}\label{subsec_orbitraising}

The first test case considered is an orbit raising scenario in which the spacecraft moves from LEO to GEO, and the final polar angle is free. For comparison, the same design parameters and boundary conditions from Taheri and Abdelkhalik are used \cite{doi:10.2514/1.58789}. The initial LEO and final GEO orbits are circular, and the boundary conditions are shown in Table \ref{tab:wob_ic_2_3_40}. The gravitational parameter for the Earth is $\mu=398601.2$ km$^3$/s$^2$. The angular rates of the spacecraft in the LEO and GEO orbits are then computed from $\dot{\theta}=\sqrt{\mu/r^3}$. The boundary conditions and parameters are converted to canonical units prior to simulation. All of the figures which follow are presented using canonical units.

The number of Fourier terms used to approximate the radial and angular states as well as the number of equally spaced discrete points (DP) along the trajectory are shown in Table \ref{tab:Design_Parameters}. The maximum and minimum thrust accelerations allowed throughout the flight are designated as $T_{a,max}$, and the time of flight weight is $\omega$.

\begin{table}[htbp]
	\fontsize{10}{10}\selectfont
    \caption{Boundary Conditions: Orbit Raising}
   \label{tab:wob_ic_2_3_40}
        \centering 
   \begin{tabular}{c | c } % Column formatting, 
      \hline 
        Parameter & Value \\ 
        \hline
        $r_i$ & 6570 km \\ 
        $\theta_i$ & 0 deg. \\ 
        $\dot{r}_i$ & 0 km/s \\
        $r_f$ & 42160 km \\ 
        $\dot{r}_f$ & 0 km/s \\
      \hline
   \end{tabular}
\end{table}

\begin{table}[htbp]
	\fontsize{10}{10}\selectfont
    \caption{Design Parameters: Orbit Raising Case 1}
   \label{tab:Design_Parameters}
        \centering 
   \begin{tabular}{c | c } % Column formatting, 
      \hline 
        Parameter & Value \\ 
        \hline
        $n_r$ & 2 \\ 
        $n_\theta$ & 3 \\ 
        $DP$ & 160 \\
        $T_{a,max}$ & $\pm$ 0.0102 DU/TU$^2$ \\
        $\omega$ & 0.5 \\
      \hline
   \end{tabular}
\end{table}

\begin{figure}[h]
  \centering
  \begin{minipage}[b]{0.45\textwidth}
    \includegraphics[width=\textwidth]{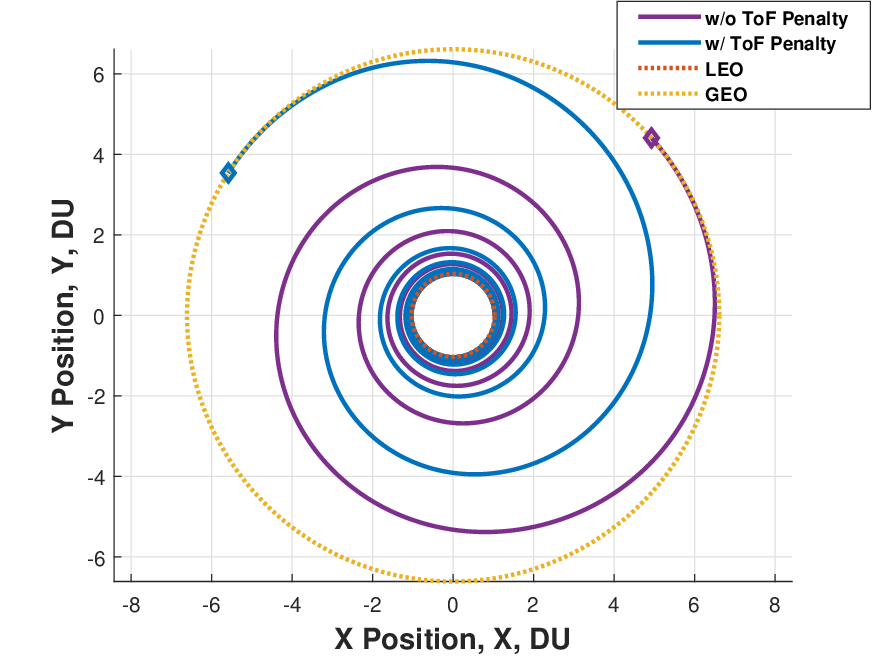}
    \caption{FFS Trajectory: Orbit Raising}
    \label{fig:trajectory}
  \end{minipage}
  \hfill
  \begin{minipage}[b]{0.45\textwidth}
    \includegraphics[width=\textwidth]{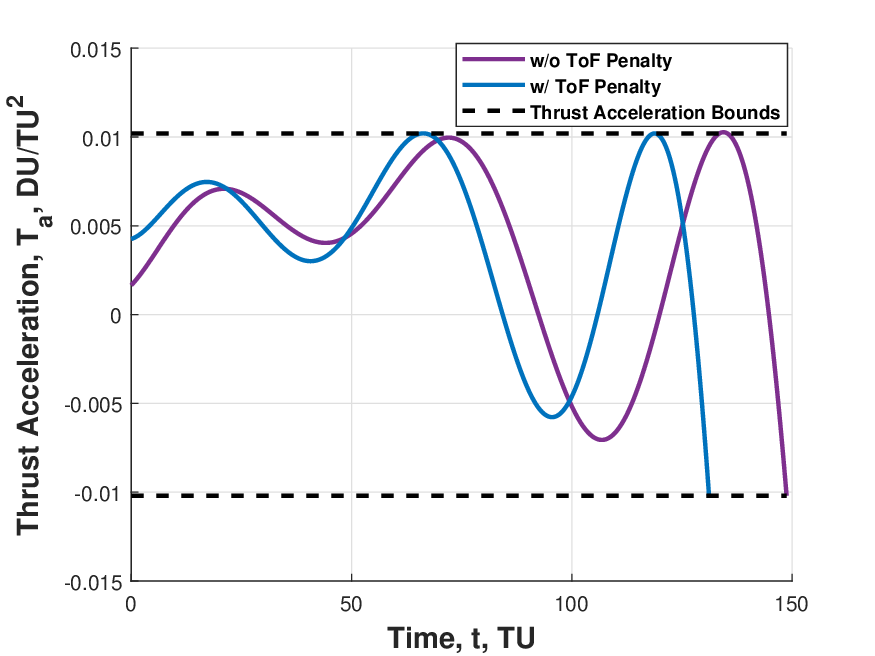}
    \caption{Thrust Acceleration: Orbit Raising}
    \label{fig:Thrust_Acceleration_Profile}
  \end{minipage}
\end{figure}

The weight for the time of flight penalty is 50\% so both the equation of motion error and time of flight objectives are weighted equally. The low-thrust FFS SB trajectory generated from penalizing the ToF is compared against previously verified results shown below \cite{doi:10.2514/1.58789}. The trajectory and thrust acceleration shown in Figure \ref{fig:trajectory} and Figure \ref{fig:Thrust_Acceleration_Profile} show the time of flight decreases as expected when the time of flight is penalized. The trajectory resulted in a ToF reduction of $\Delta ToF = -3.95$ hours - an 11.84\% reduction. However, the equation of motion error increased which indicates the FFS approximation does not satisfy the equations of motion as well.  

\begin{figure}[h]
  \centering
  \begin{minipage}[b]{0.45\textwidth}
    \includegraphics[width=\textwidth]{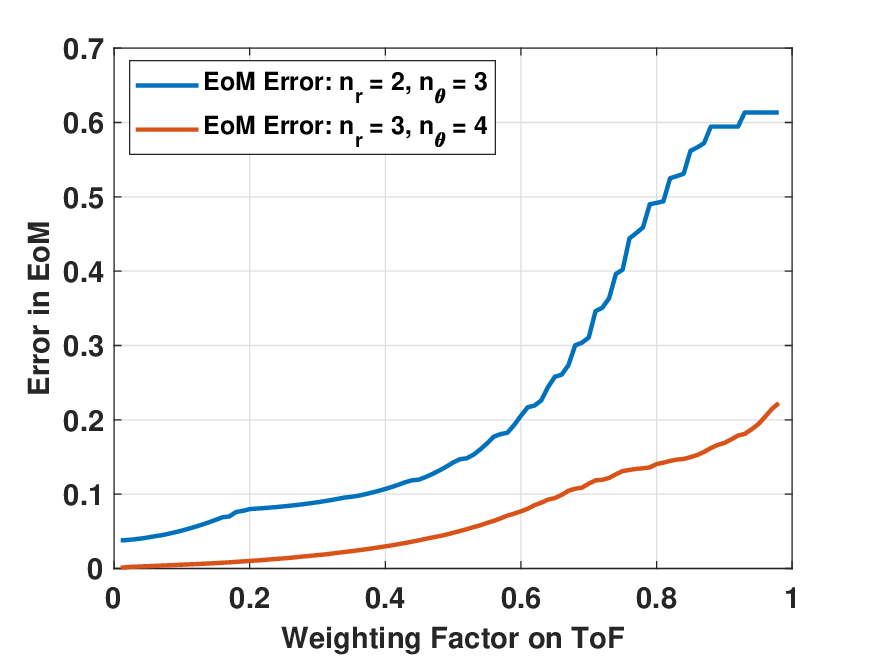}
    \caption{Optimal EoM vs. Weight}
    \label{fig:EOM_VS_WF}
  \end{minipage}
  \hfill
  \begin{minipage}[b]{0.45\textwidth}
    \includegraphics[width=\textwidth]{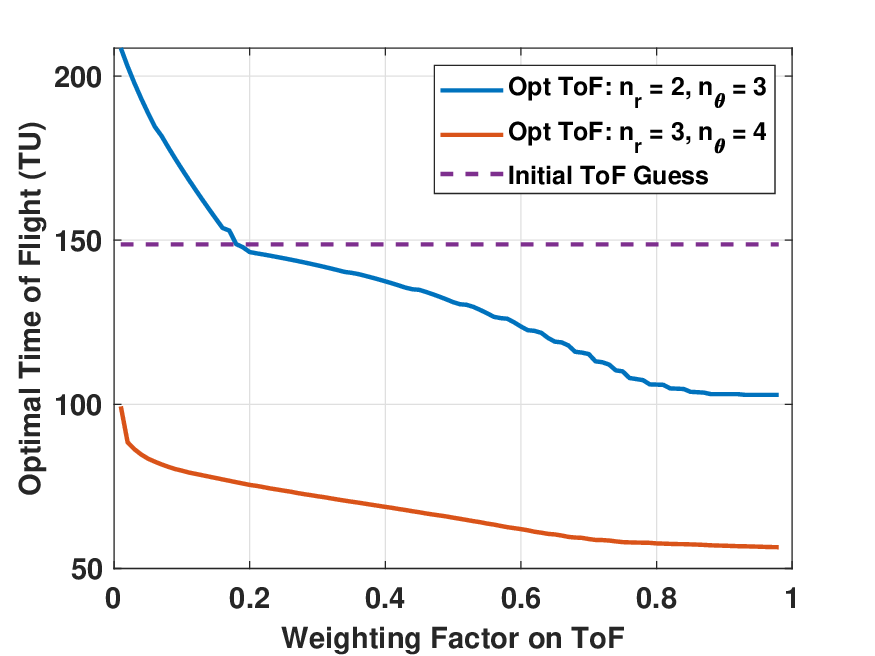}
    \caption{Optimal ToF vs Weight}
    \label{fig:TOF_VS_WF}
  \end{minipage}
\end{figure}

The reduction in time of flight (at the cost of additional equation of motion error) is an improvement, but the time of flight weights which produce feasible trajectories must still be determined. The results above are extended to show the trade off between the equation of motion error and the time of flight objectives by plotting both as functions of the time of flight weight in Figure \ref{fig:EOM_VS_WF} and Figure \ref{fig:TOF_VS_WF}. The figures include plots of the objectives for trajectories using FFS SB parameters $n_r=2$ and $n{\theta}=3$ as well as $n_r=3$ and $n_{\theta}=4$. The trajectories were generated for weights between 0.01 and 0.99 in steps of 0.01. 

Figure \ref{fig:EOM_VS_WF} shows that error increases as the weight increases, which is expected. At low weights the error in the EoM is very small but grows continually as the normalized ToF becomes increasingly more important in the objective function. Another important observation is that increasing the number of Fourier coefficients to approximate the trajectory decreases the total error which is shown in Figure \ref{fig:EOM_VS_WF}, and reduces the total time of flight which is shown in Figure \ref{fig:TOF_VS_WF}. These plots were used to determine the trade-offs in the objective function for the purpose of selecting better time of flight weights which satisfy the multi-objective optimization.

\subsubsection{Trajectory Feasibility}

The FFS trajectory results are compared against results generated from numerically integrating the equations of motion. A weight of of 0.01 was selected for the following test cases. It is completely viable to choose a larger weight, however, the drawback is that the error in the equations of motion will increase. This will result in larger position errors in the spacecraft trajectory FFS approximation.

The numerical integration was implemented using the built in MATLAB function \textit{ode45}. This function is used to solve the system of differential equations for the spacecraft shown in Eqs. \eqref{eq:r_eom} and \eqref{eq:theta_eom}. The simulation uses an open loop control approach in which the FFS SB thrust acceleration profile is interpolated at each time step. The design parameters used to obtain the following results are shown in Table \ref{tab:Design_Parameters2}.

\begin{table}[htbp]
	\fontsize{10}{10}\selectfont
    \caption{Design Parameters: Orbit Raising Case 2}
   \label{tab:Design_Parameters2}
        \centering 
   \begin{tabular}{c | c } % Column formatting, 
      \hline 
        Parameter & Value \\ 
        \hline
        $n_r$ & 3 \\ 
        $n_\theta$ & 4 \\ 
        $DP$ & 160 \\
        $T_{a,max}$ & $\pm$ 0.0102 DU/TU$^2$ \\
        $\omega$ & 0.01 \\
      \hline
   \end{tabular}
\end{table}

\begin{figure}[h]
  \centering
  \begin{minipage}[b]{0.45\textwidth}
    \includegraphics[width=\textwidth]{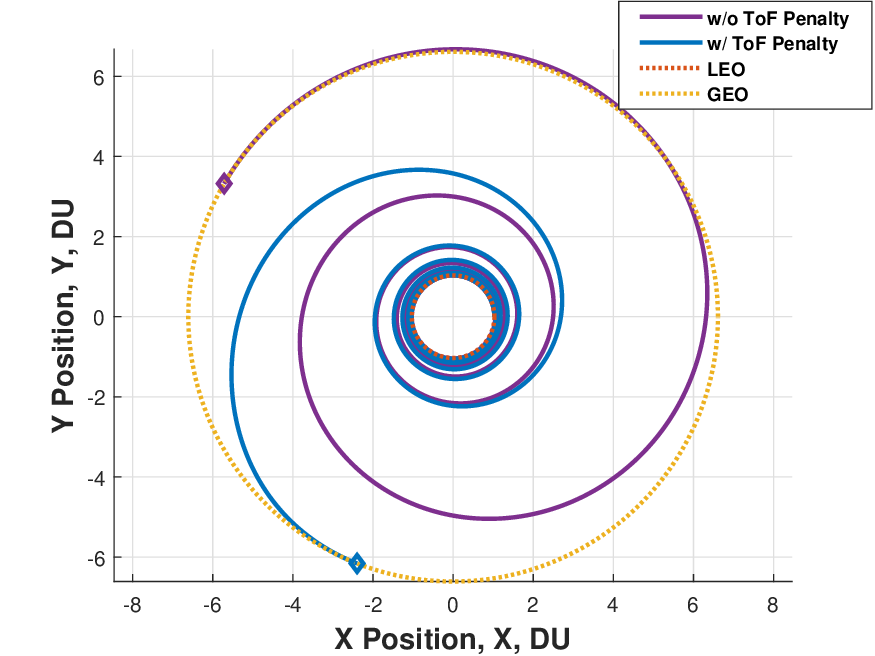}
    \caption{FFS Trajectories: Orbit Raising}
    \label{fig:wob_traj_3_4_160_1}
  \end{minipage}
  \hfill
  \begin{minipage}[b]{0.45\textwidth}
    \includegraphics[width=\textwidth]{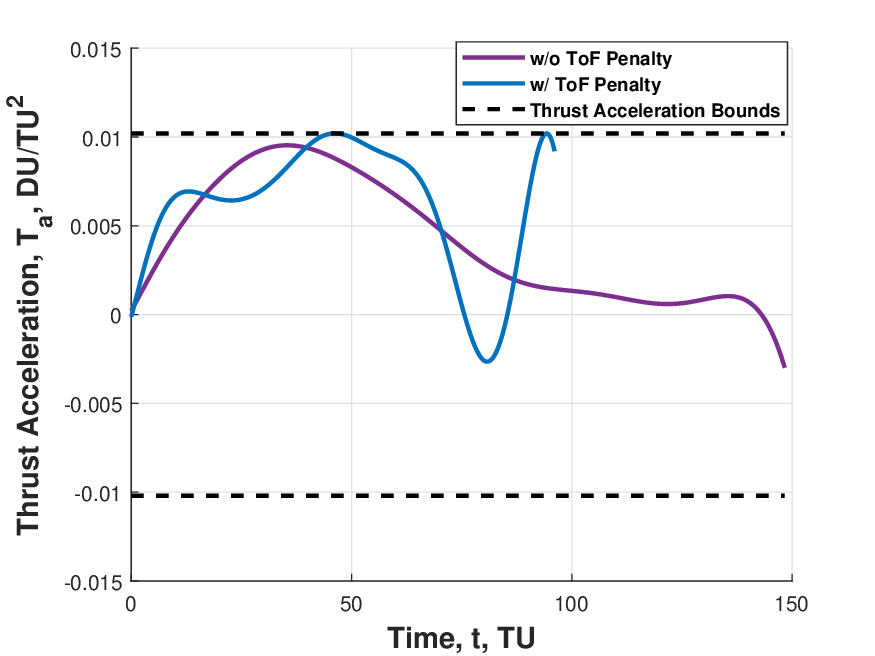}
    \caption{\small{Thrust Acceleration: Orbit Raising}}
    \label{fig:wob_ta_3_4_160_1}
  \end{minipage}
\end{figure}

In Figure \ref{fig:wob_traj_3_4_160_1} two FFS trajectories are shown. For one trajectory the time of flight weight is $\omega=0$, meaning the only term contributing to the objective function is the term containing the error in the EoM, $\textbf{F}^T \textbf{F}$. The other trajectory has a time of flight weight of $\omega=0.01$. Comparing the two trajectories it can be seen that the spacecraft spends much more time closer to the Earth on all but the final revolution of the penalized trajectory. On this final revolution, the radial position increases quickly to the desired radius rather than slowly growing over multiple revolutions as it does in the non-penalized trajectory.

\begin{figure}[h]
  \centering
  \begin{minipage}[b]{0.5\textwidth}
    \includegraphics[width=\textwidth]{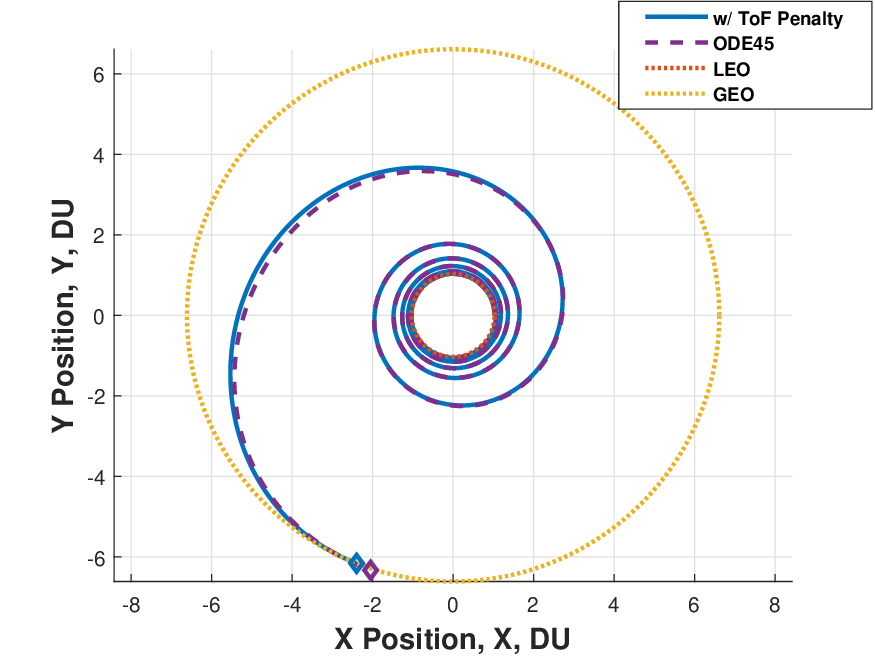}
    \caption{\small{Trajectory Comparison: Orbit Raising}}
    \label{fig:wob_traj_3_4_160_1_ode45}
  \end{minipage}
  \hfill
\end{figure}

The thrust acceleration curves for these two trajectories are shown in Figure \ref{fig:wob_ta_3_4_160_1} and it can be seen that the penalized trajectory takes much less time to reach the desired orbit. The equation of motion error between the FFS approximation and the \textit{ode45} trajectory can be seen in Figure \ref{fig:wob_traj_3_4_160_1_ode45}. The trajectories and final polar angles are slightly different, however, the trajectories are still nearly a match.

%start here for new contributions for orbit raising

\begin{figure}[h]
  \centering
  \begin{minipage}[b]{0.51\textwidth}
    \includegraphics[width=\textwidth]{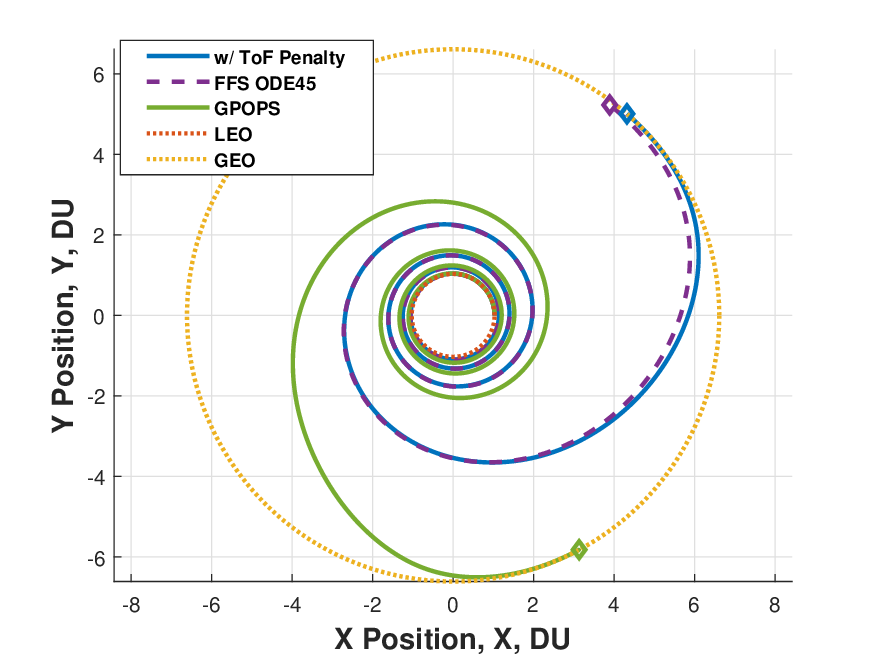}
    \caption{FFS and GPOPS Trajectories: Orbit Raising}
    \label{fig:wob_traj_5_7_40_1}
  \end{minipage}
  \hfill
  \begin{minipage}[b]{0.45\textwidth}
    \includegraphics[width=\textwidth]{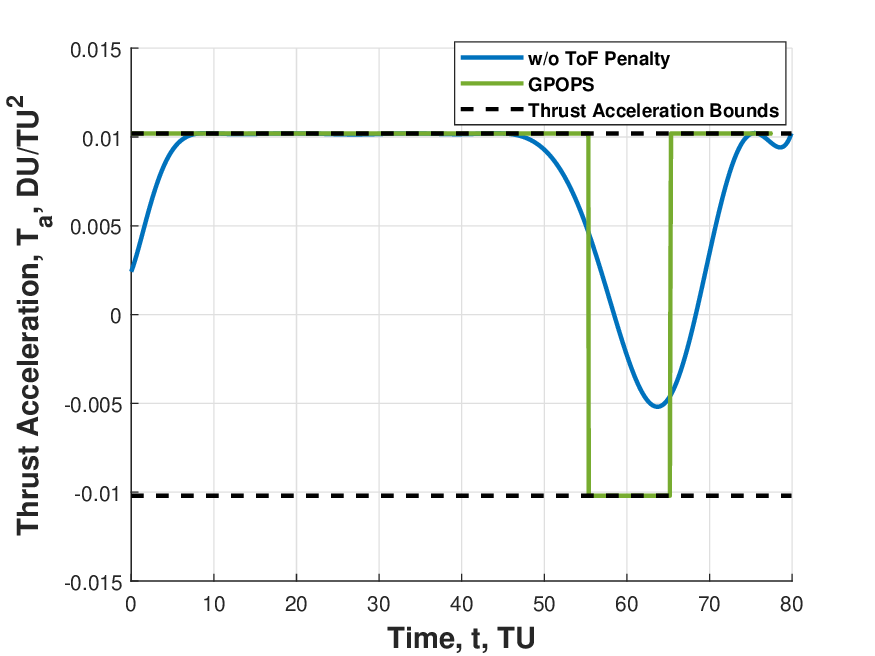}
    \caption{Thrust Acceleration: Orbit Raising}
    \label{fig:wob_ta_5_7_40_1}
  \end{minipage}
\end{figure}

One final orbit raising trajectory was generated with the parameters shown in Table \ref{tab:Design_Parameters3}. This test case uses more $n_r$ and $n_\theta$ terms. In addition, the number of discretization points had to be decreased from 160 points to 40 points in order to generate Fourier coefficients which produced a feasible trajectory. Results are shown in figures \ref{fig:wob_traj_5_7_40_1} and \ref{fig:wob_ta_5_7_40_1}. This test case was specifically included to demonstrate how close the trajectories generated using FFS can be to the optimal solutions generated by GPOPS.

\begin{table}[htbp!]
	\fontsize{10}{10}\selectfont
    \caption{Design Parameters: Orbit Raising Case 3}
   \label{tab:Design_Parameters3}
        \centering 
   \begin{tabular}{c | c } % Column formatting, 
      \hline 
        Parameter & Value \\ 
        \hline
        $n_r$ & 5 \\ 
        $n_\theta$ & 7 \\ 
        $DP$ & 40 \\
        $T_{a,max}$ & $\pm$ 0.0102 DU/TU$^2$ \\
        $\omega$ & 0.01 \\
      \hline
   \end{tabular}
\end{table}

% ---------- GPOPS ----------
% Youngro take here
% ---------------------------

In addition to the FFS SB and associated \textit{ode45} trajectories, figures \ref{fig:wob_traj_5_7_40_1} and \ref{fig:wob_ta_5_7_40_1} show the optimal trajectory and thrust acceleration profile generated using GPOPS.
The optimal flight time is 17.38 hours, which corresponds to a 47.87 \% reduction compared to the FFS trajectory with no ToF weight for case 3. This is very close to the FFS trajectory time of flight, $ToF=17.91$ hours, for the FFS parameters in Table \ref{tab:ToForbitraising}. Even though the optimal time of flight and the FFS time of flight are very close, the trajectories are still quite different due to variance in their thrust acceleration profiles. Given these differences, the optimal trajectory orbits almost 4 times while the FFS trajectory orbits 5 times.
The most distinctive feature of the optimal thrust acceleration solution is that it switches its value from one extreme to the other, which is referred to as a bang-bang solution.

All results for the orbit raising cases presented above are shown in Table \ref{tab:ToForbitraising}. The GPOPS optimal ToF is compared to the non-penalized ToF result generated in case 3, $ToF=33.34$ hours.

\begin{table}[htbp!]
	\fontsize{10}{10}\selectfont
    \caption{Time of Flight Results: Orbit Raising}
   \label{tab:ToForbitraising}
        \centering 
   \begin{tabular}{c | c | c | c} % Column formatting, 
      \hline 
        Parameter & Case 2 & Case 3 & GPOPS\\
        \hline
        $n_r$ & 3 & 5 & - \\
        $n_{\theta}$ & 4 & 7 & - \\
        $DP$ & 160 & 40 & - \\
        $\omega$ & 0.01 & 0.01 & - \\ 
        \hline
        ToF w/o Penalty (hr)& 33.24 & 33.34 & - \\
        ToF w/ Penalty (hr)& 21.53 & 17.91 & - \\
        TOF GPOPS (hr) & - & - & 17.38 \\
        ToF Reduction (\%) & 35.23 & 46.29 & 47.87\\
      \hline
   \end{tabular}
\end{table}

\subsection{Rendezvous}\label{subsec_rendezvous}

The same process was used to compare the trajectories obtained for the rendezvous case. In addition to the boundary conditions given in Table \ref{tab:ic_3_4_160_ren}, the same design parameters from the orbit raising test case were used to produce the results for this test case shown in Table \ref{tab:Design_Parameters4}. It is worth noting that identical design parameters do not necessarily result in a good spacecraft trajectory for each maneuver, however, in this case they did. In general, different scenarios and boundary conditions will require unique FFS design parameters.

\begin{table}[htbp!]
	\fontsize{10}{10}\selectfont
    \caption{Boundary Conditions: Rendezvous}
   \label{tab:ic_3_4_160_ren}
        \centering 
   \begin{tabular}{c | c } % Column formatting, 
      \hline 
        Parameter & Value \\ 
        \hline
        $r_i$ & 6570 km \\ 
        $\theta$ & 0 deg. \\ 
        $\dot{r}_i$ & 0 km/s \\
        $r_f$ & 42160 km \\ 
        $\dot{r}_f$ & 0 km/s \\
      \hline
   \end{tabular}
\end{table}

\begin{table}[htbp!]
	\fontsize{10}{10}\selectfont
    \caption{Design Parameters: Rendezvous Case 4}
   \label{tab:Design_Parameters4}
        \centering 
   \begin{tabular}{c | c } % Column formatting, 
      \hline 
        Parameter & Value \\ 
        \hline
        $n_r$ & 3 \\
        $n_\theta$ & 4 \\
        $DP$ & 160 \\
        $T_{a,max}$ & $\pm$ 0.0102 DU/TU$^2$ \\
        $\omega$ & 0.01 \\
      \hline
   \end{tabular} 
\end{table}

\begin{figure}[h]
  \centering
  \begin{minipage}[b]{0.45\textwidth}
    \includegraphics[width=\textwidth]{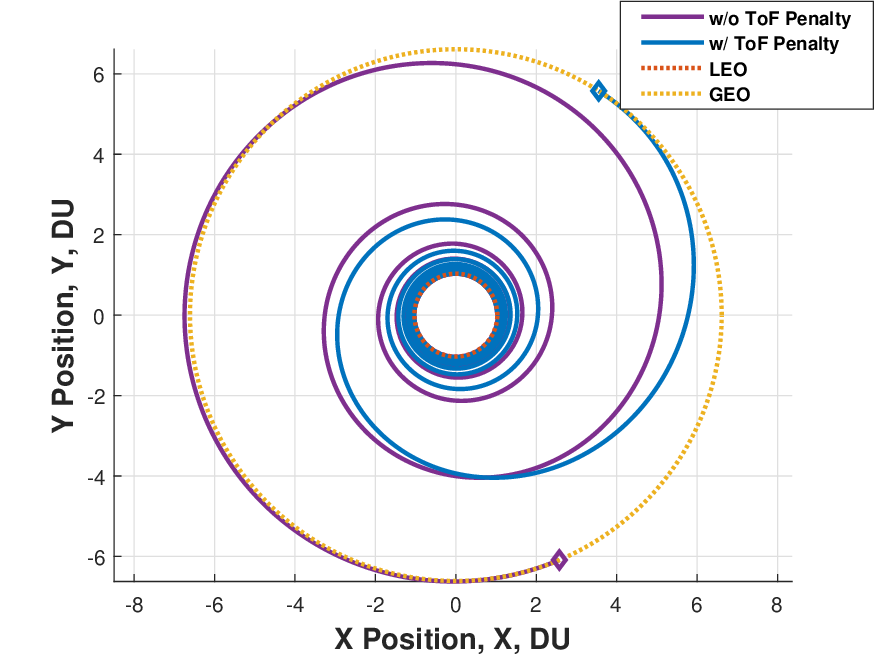}
    \caption{FFS Trajectories: Rendezvous}
    \label{fig:wr_traj_3_4_160_1}
  \end{minipage}
  \hfill
  \begin{minipage}[b]{0.45\textwidth}
    \includegraphics[width=\textwidth]{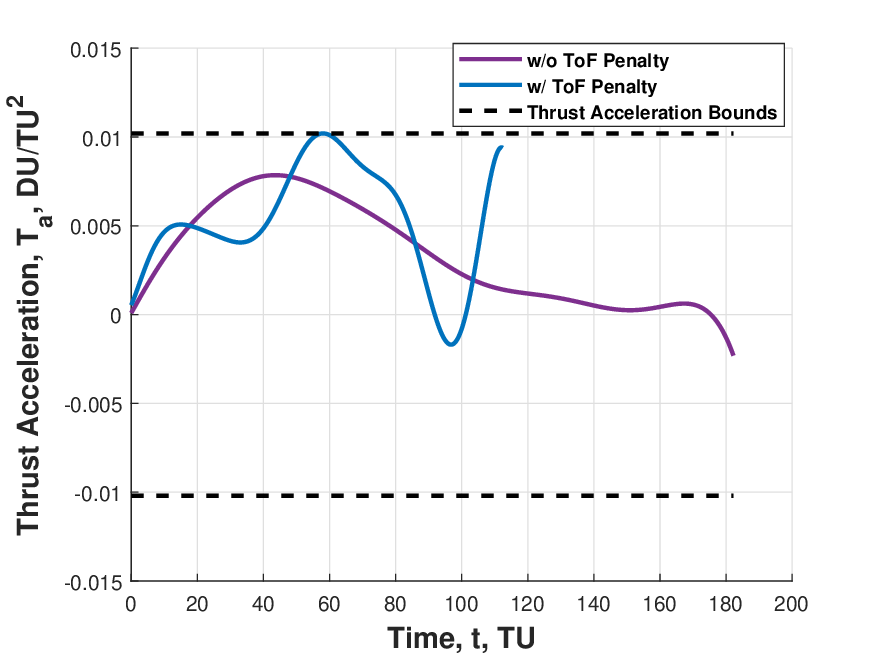}
    \caption{Thrust Acceleration: Rendezvous}
    \label{fig:wr_ta_3_4_160_1}
  \end{minipage}
\end{figure}

\begin{figure}[h!]
  \centering
  \begin{minipage}[b]{0.45\textwidth}
    \includegraphics[width=\textwidth]{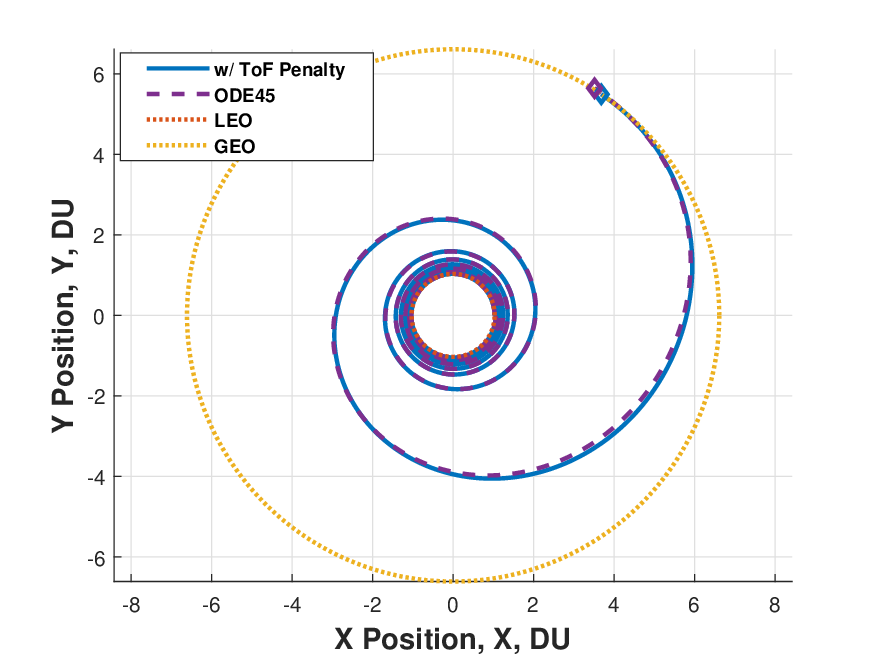}
    \caption{Trajectory Comparison: Rendezvous}
    \label{fig:wr_traj_3_4_160_1_ode45}
  \end{minipage}
  \hfill
\end{figure}

Two FFS trajectories, with and without penalties on the ToF, were generated and are shown in Figure \ref{fig:wr_traj_3_4_160_1}. The thrust profiles for these two trajectories are shown in Figure \ref{fig:wr_ta_3_4_160_1} and it is clear that the penalized trajectory completes the rendezvous in much shorter time. A larger thrust acceleration must be applied to rendezvous in a shorter time period. Figure \ref{fig:wr_traj_3_4_160_1_ode45} shows the FFS trajectory compares well to the trajectory obtained via numerically integrating the thrust profile from the FFS approximation.

\subsubsection{Trajectory Feasibility}
% ---------- GPOPS ----------
% Youngro take here
% ---------------------------

\begin{figure}[h]
  \centering
  \begin{minipage}[b]{0.45\textwidth}
    \includegraphics[width=\textwidth]{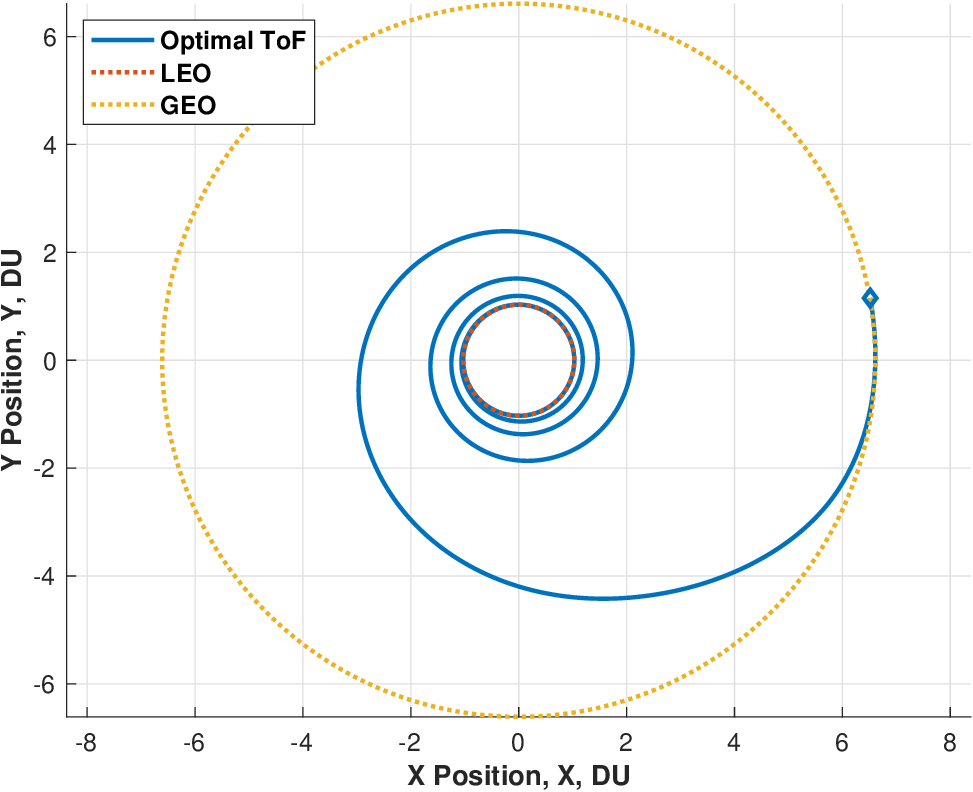}
    \caption{GPOPS Trajectory: Rendezvous}
    \label{fig:gpops_rendezvous_trajectory}
  \end{minipage}
  \hfill
  \begin{minipage}[b]{0.50\textwidth}
    \includegraphics[width=\textwidth]{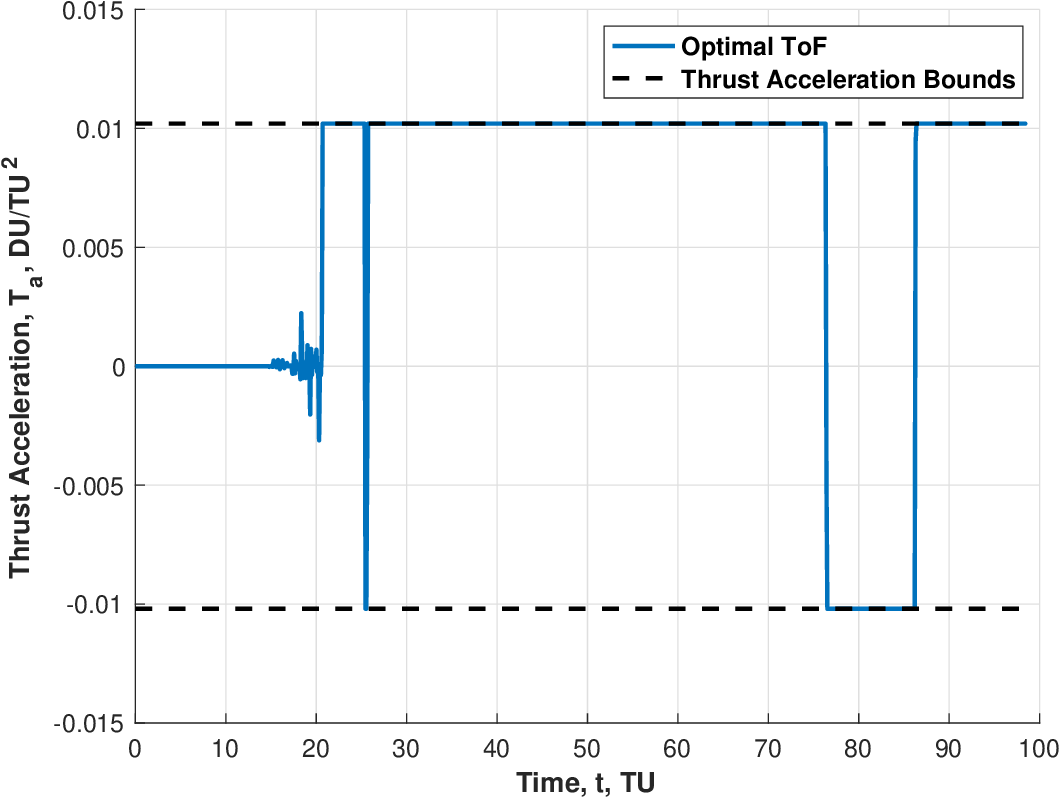}
    \caption{GPOPS Thrust Acceleration: Rendezvous}
    \label{fig:gpops_rendezvous_thrust}
  \end{minipage}
\end{figure}

The optimal trajectory and thrust acceleration profile that achieves the minimum time rendezvous are shown in Figures \ref{fig:gpops_rendezvous_trajectory} and \ref{fig:gpops_rendezvous_thrust}.
Contrary to the previous case, the optimal thrust acceleration for the rendezvous case has a bang-off-bang structure.
The solution stays near zero until 15 TU, oscillations occur around 20 TU, and the thrust profile touches the lower bound for an instant around 25 TU. 
This oscillation may be attributed to the fact that the final polar angle is updated every time a solver inside GPOPS is called.
As a result, an optimal time of flight of 22.08 hours is obtained which corresponds to a 45.94 \% reduction compared to the non-weighted FFS trajectory time of flight.

All results for the rendezvous cases presented above are shown in table \ref{tab:ToFrendezvous}. The GPOPS optimal ToF is compared to the non-penalized ToF result that is generated in case 4, $ToF=40.84$ hours.

\begin{table}[h!]
	\fontsize{10}{10}\selectfont
    \caption{Time of Flight Results: Rendezvous}
   \label{tab:ToFrendezvous}
        \centering 
   \begin{tabular}{c | c | c} % Column formatting, 
      \hline 
        Parameter & Case 4 & GPOPS\\
        \hline
        $n_r$ & 3 & - \\
        $n_{\theta}$ & 4 & - \\
        $DP$ & 160 & - \\
        $\omega$ & 0.01 & - \\ 
        \hline
        ToF w/o Penalty (hr)& 40.84 & - \\
        ToF w/ Penalty (hr)& 25.19 & - \\
        ToF GPOPS (hr) & - & 22.08 \\
        ToF Reduction (\%) & 38.32 & 45.94\\
      \hline
   \end{tabular}
\end{table}

\section{Conclusion}\label{conclusion}
This paper built on the FFS SB method developed to approximate low-thrust spacecraft trajectories by adding time of flight to the objective function\cite{doi:10.2514/1.58789}. A weight is applied to the terms in the objective function to balance the reduction in the time of flight with a trajectory that satisfies the equations of motion. The time of flight penalty was tested on orbit raising and rendezvous maneuvers. In both cases the new objective function generated trajectories which decreased the ToF. The reduced time of flight trajectories were validated against numerically integrated results using the non-linear equations of motion. Additionally, the results from FFS were compared to the results from GPOPS to see how well the FFS trajectories and thrust accelerations matched the optimal solution. The results showed that the FFS obtained solutions that were close to the optimal solutions generated by GPOPS for the test cases investigated.

%\bibliographystyle{AAS_publication}   % Number the references.
%\bibliography{references}   % Use references.bib to resolve the labels.

% \bibliographystyle{unsrt}
% \bibliography{references}

\begin{thebibliography}{}

\end{thebibliography}


\begin{thebibliography}{10}

\bibitem{doi:10.2514/2.4231}
John~T. Betts.
\newblock Survey of numerical methods for trajectory optimization.
\newblock {\em Journal of Guidance, Control, and Dynamics}, 21(2):193--207,
  1998.

\bibitem{doi:10.2514/1.58789}
Ehsan Taheri and Ossama Abdelkhalik.
\newblock Shape based approximation of constrained low-thrust space
  trajectories using fourier series.
\newblock {\em Journal of Spacecraft and Rockets}, 49(3):535--546, 2012.

\bibitem{doi:10.2514/1.G000878}
Ehsan Taheri and Ossama Abdelkhalik.
\newblock Fast initial trajectory design for low-thrust restricted-three-body
  problems.
\newblock {\em Journal of Guidance, Control, and Dynamics}, 38(11):2146--2160,
  2015.

\bibitem{doi:10.1016/j.asr.2015.11.034}
Ehsan Taheri and Ossama Abdelkhalik.
\newblock Initial three-dimensional low-thrust trajectory design.
\newblock {\em Advances in Space Research}, 57(3):889--903, 2016.

\bibitem{doi:10.2514/1.36848}
Bradley~J. Wall and Bruce~A. Conway.
\newblock Shape-based approach to low-thrust rendezvous trajectory design.
\newblock {\em Journal of Guidance, Control, and Dynamics}, 32(1):95--101,
  2009.

\bibitem{doi:10.2514/8.7261}
D.~J. BENNEY.
\newblock Escape from a circular orbit using tangential thrust.
\newblock {\em Journal of Jet Propulsion}, 28(3):167--169, 1958.

\bibitem{doi:10.2514/3.56583}
Frederick~W. Boltz.
\newblock Orbital motion under continuous tangential thrust.
\newblock {\em Journal of Guidance, Control, and Dynamics}, 15(6):1503--1507,
  1992.

\bibitem{doi:10.2514/1.13095}
Anastassios~E. Petropoulos and James~M. Longuski.
\newblock Shape-based algorithm for the automated design of low-thrust, gravity
  assist trajectories.
\newblock {\em Journal of Spacecraft and Rockets}, 41(5):787--796, 2004.

\bibitem{petrothesis}
A.E. Petropoulos.
\newblock A shape-based approach to automated, low-thrust, gravity-assist
  trajectory design.
\newblock {\em Ph.D. Thesis}, 2001.

\bibitem{patterson2014gpops}
Michael~A Patterson and Anil~V Rao.
\newblock Gpops-ii: A matlab software for solving multiple-phase optimal
  control problems using hp-adaptive gaussian quadrature collocation methods
  and sparse nonlinear programming.
\newblock {\em ACM Transactions on Mathematical Software (TOMS)}, 41(1):1--37,
  2014.

\end{thebibliography}

\end{document}